\documentclass[graybox]{svmult}

\usepackage{type1cm}        
\usepackage{makeidx}         
\usepackage{graphicx}        
\usepackage{multicol}        
\usepackage[bottom]{footmisc}

\usepackage{newtxtext}       %
\usepackage[varvw]{newtxmath}       
\usepackage{epstopdf}
\usepackage{caption}
\usepackage{subcaption}

\newcommand{\U}{\mathbf{U}}

\def\bcdot{{*}}

\newcommand{\ha}{\frac{1}{2}}

\usepackage[normalem]{ulem}

\def\be{\begin{equation}}
\def\ee{\end{equation}}
\newcommand{\ei}[0]{\end{itemize}}
\newcommand{\beann}[0]{\begin{eqnarray*}}
\newcommand{\eeann}[0]{\end{eqnarray*}}
\def\bea{\begin{eqnarray}}
\def\eea{\end{eqnarray}}
\def\ba{\begin{array}{l}\displaystyle}
\def\ea{\end{array}}

\newcommand{\aposteriori}{\textit{a posteriori} }

\newcommand{\ema}[1]{{\color{black} #1}}
\newcommand{\raph}[1]{{\color{black} #1}}
\newcommand{\gio}[1]{{\color{black} #1}}
\newcommand{\add}[1]{{\color{black} #1}}

\makeindex             

\begin{document}

\title{CAT-MOOD methods for conservation laws in one space dimension}
\author{R. Loubere, E. Macca, C. Pares and G. Russo}
\institute{Emanuele Macca, corresponding author, \at Universit{y} of Catania, \email{emanuele.macca@unict.it}
\and Rapha{\"e}l Loub{\`e}re \at  Universit{y} of Bordeaux, CNRS, \email{raphael.loubere@math.u-bordeaux.fr} 
\and Carlos Par{\'e}s \at  University of M{\'a}laga, \email{pares@uma.es}
\and Giovanni Russo \at Universit{y} di Catania, \email{russo@dmi.unict.it}
}

%
\maketitle

\abstract{In this paper we blend high-order Compact Approximate Taylor (CAT) numerical methods with the \textit{a posteriori} Multi-dimensional Optimal Order Detection (MOOD) paradigm to solve hyperbolic systems of conservation laws. The resulting methods 
are highly accurate for smooth solutions, essentially non-oscillatory for discontinuous ones, and almost fail-safe positivity preserving. 
Some numerical results for scalar conservation laws and systems are presented to show the appropriate behavior of CAT-MOOD methods.
}


\section{Introduction} \label{sec:introduction}
Lax-Wendroff methods  for linear systems of conservation laws are based on Taylor expansions in time in which the time derivatives are transformed into spatial derivatives using the governing equations  \cite{LeVeque2007book,Toro2009}. The spatial derivatives are then discretised by means of centered high-order differentiation formulas. 

 \add{One of the} 
difficulties to extend Lax-Wendroff methods to nonlinear problems come from the transformation of time derivatives into spatial derivatives through the Cauchy-Kovalesky (CK) procedure: this approach may indeed be impractical from the computational point of view. The Lax-Wendroff Approximate Taylor (LAT) methods introduced in  \cite{ZBM2017}  circumvent the CK procedure by computing  time derivatives in a recursive way using high-order centered differentiation formulas combined with  Taylor expansions in time. Compact Approximated Taylor methods (CAT) introduced in \cite{Carrillo-Pares}  follow a similar strategy.  These methods are compact in the sense that the length of the stencils is minimal: $(2P +1)$-point stencils are used to get order $2P$ compared to $4P +1$-point stencils in LAT methods. They are also  $L^2$ linearly-stable.

\add{A second difficulty comes from the treatment of shocks and discontinuities that usually arise in quasilinear systems of conservation laws.} 
In order to avoid the spurious oscillations that Lax-Wendroff-type methods produce in presence of discontinuities or high gradients, CAT methods were combined in \cite{CPZMR2020,TesiPhD,Macca-Pares} with an \textit{a priori} order adaptive procedure.
\raph{To do so} a family of smoothness indicators were used to automatically reduce the order of the method 
\raph{in the vicinity of}
discontinuities. 
\raph{These limited CAT methods are called ACAT.}

The goal of this paper is to combine \raph{1D} CAT methods with the \textit{a posteriori} Multi-dimensional Optimal Order Detection (MOOD) paradigm introduced in \cite{CDL1,CDL0_FVCA}. \ema{This technique is expected to produce non-oscillatory high-order methods with an appropriate detection of discontinuities} The \raph{resulting} methods will be applied to scalar conservation laws and the 1D Euler equations of gas-dynamics.

The rest of this paper is organized as follows: 
the \add{next} section introduces the system of equations we plan to solve. 
In section three, CAT methods are recalled. 
The fourth section describes how to blend CAT schemes with MOOD.
\add{Numerical results are reported in the fifth section, illustrating the behavior} of the sixth-order CAT-MOOD method.
Conclusions and perspectives are finally drawn.

%
\section{Governing equations} \label{sec:equations}
We consider 1D hyperbolic systems of conservation laws of the form
\bea \label{eq:systemPDEs_bis}
    \partial_t \U + \partial_x \mathbf{F}(\U) = \mathbf{0},
\eea
where $t\in \mathbb{R}^+$ represents the time variable, $x\in \mathbb{R}$ the space variable and $\U=\U(x,t)
\in \mathbb{R}^M$ is the vector of conserved variables while 
$\mathbf{F}(\U(x,t)) \in \mathbb{R}^M$ is the flux vector. More precisely, we focus on the 1D Euler equations of gas dynamics in which $M = 3$,  $\U=(\rho, \rho u,  \rho e)^t$ with $\rho$ the density, $u$ the velocity and 
$e=\varepsilon + \frac12 \raph{u^2}$ the total energy \add{per unit mass}, being $\varepsilon$ the specific \raph{internal energy}.
The flux is given by 
$\mathbf{F}(\U)= \left(
\rho u, \rho u^2 +p, (\rho e +p )u\right)^t$.
The system is closed with the \raph{perfect gas} equation of state:
$ \raph{p(\rho,\varepsilon) = (\gamma-1)\rho \varepsilon}. $
The equations in (\ref{eq:systemPDEs_bis}) represent the conservation of mass, momentum and total energy.
An entropy inequality has to be 
\raph{supplemented to deal with discontinuous solutions}.
This system is hyperbolic with eigenvalues $\lambda^-=u - c$, $\lambda^0=u,$  $\lambda^+=u + c$ where $c = \sqrt{\gamma p/\rho}$ \raph{is the sound speed}.
The set of admissible states, i.e. states that are consistent with physics, is
\bea \label{eq:physical_states}
    \mathcal{A} = \left\{ \U\in \mathbb{R}^M, \; \text{such that} \; \rho>0, \; p>0 \right\}.
\eea
Together with this system, two scalar conservation laws
\raph{of the form}
\begin{equation}
    \label{sec:CAT_gov_equ}
    u_t + \partial_x f(u) = 0,
\end{equation} 
will be considered 
to test the methods: \raph{first} the linear advection equation corresponding to $f(u)=b u$ with $b\in \mathbb{R}$, and, \raph{secondly}, Burgers' equation corresponding to $f(u)=u^2/2$.  
\raph{These scalar conservation laws verify the maximum principle, that is
\begin{eqnarray} \label{eq:physical_states2}
    \mathcal{A} & = & \left\{ u\in \mathbb{R}, \; \text{s.t.} \; \underline{u}_0 \leq u(x,t) \leq \overline{u}_0 \right\}, 
    \nonumber\\
    \underline{u}_0 & = & \min_x( u(x,t_0)), \; \overline{u}_0=\max_x( u(x,t_0)), 
\end{eqnarray}
where $t_0$ is the initial time, and $u(x,t_0)$ the initial condition.
}

\section{Compact Approximate Taylor (CAT) schemes} \label{sec:CAT} 
In this section we recall the formulation of CAT methods. In order to simplify the notation, the methods are introduced for scalar conservation laws \eqref{sec:CAT_gov_equ}, but the expression for systems \eqref{eq:systemPDEs_bis} is similar.

We consider uniform meshes in 1D: the space domain is split into computational cells  $\omega_i=[x_{i-1/2},x_{i+1/2}]$ of constant  \add{width} $\Delta x$. $x_i=\ha(x_{i+1/2}+x_{i-1/2})$ represents the center of the $i$-th cell.
Although in practice the time step depends on the CFL condition and thus it is not constant, for the sake of simplicity in the presentation of the methods it will be assumed that the time interval $[0,T]$ is split into sub-intervals $[t^n,t^{n+1}]$ of constant length $\Delta t$.

The $2P$-CAT method is written in conservative form as
\bea \label{eq:general_LW}
 u_i^{n+1} = u_i^n + \frac{\Delta t}{\Delta x}\left( F^P_{i-1/2}- F^P_{i+1/2}\right).
\eea
To compute the numerical flux $F^P_{i+1/2}$ only the approximations in the $2P$-point stencil
\bea \label{eq:stencil2P} 
\mathcal{S}_{i+\ha}^P = \left\{ x_{i-P+1},\ldots, x_{i+P} \right\}
\eea
are used, which ensures that  $(2P + 1)$ points are 
\raph{employed}
to update the numerical solution. 
The following formulas of numerical differentiation are used to compute the numerical fluxes:  given two positive integers $P$, $k$, an index $i$, and a real number $q$, we consider the interpolatory formula that approximates the $k$-th derivative of a function $f$  at the point $x_i + q \Delta x$ using its values at the $2P$  points $x_{i-P+1}, \dots, x_{i+P}$:
\begin{equation}\label{upwF}
f^{(k)}(x_i + q \Delta x) \approx  \mathcal{A}^{k,j}_{P} (f, \Delta x) = \frac{1}{\Delta x^k} \sum_{j = -P + 1}^P \gamma^{k,q}_{P,j} f(x_{i+j}).
\end{equation}
\raph{Notice that the case} $k = 0$ corresponds to Lagrange interpolation.  
When the formulas are applied to approximate the partial derivatives of a function $f(x,t)$ from some approximations $f_i^n \approx f(x_i, t_n)$, the symbol $\bcdot$ will be used to indicate to which variable (space or time) the differentiation is applied. For instance:
\begin{eqnarray*}
& &  \partial^k_x f(x_{i}, t_n)  \approx 
\mathcal{A}_{P}^{k,0} (f_{\bcdot}^n,\Delta x\Bigr)
=  \frac{1}{\Delta x^k} \sum_{j=-P+1}^{P} \gamma^{k,0}_{P,l} f_{i+j}^n, \\
& &  \partial^k_t f(x_{i}, t_n)  \approx  \mathcal{A}_{P}^{k,0}(f_i^{\bcdot},\Delta t)
 = \frac{1}{\Delta t^k} \sum_{r=-P+1}^{P} \gamma^{k,0}_{P,r} f_i^{n + r}.
 \end{eqnarray*}

Using this notation, the expression of the numerical flux is as follows:
\begin{equation}\label{cat2}
F^P_{i+1/2}  = \sum_{k=1}^{m} \frac{\Delta t^{k-1}}{k!}\mathcal{A}^{0, 1/2}_{P}(f_{i,\bcdot}^{(k-1)}, \Delta x),
\end{equation}
where 
\begin{equation}\label{f(k-1)_ij}
f^{(k-1)}_{i,j}  \approx \partial_t^{k-1}f(u)(x_{i+j}, t^n), \quad j=-P+1,\dots, P
\end{equation}
are \textit{local} approximations of the time derivatives of the flux. 
By \textit{local} we mean that these approximations depend on the stencil, i.e.
\raph{for two different stencils such that 
$i_1 + j_1 = i_2 + j_2$ then $f^{(k-1)}_{i_1,j_1}$ is not necessarily equal to $f^{(k-1)}_{i_2,j_2}$.}




Since the exact solution satisfies
$$
\partial_t^k u = - \partial_t^{k-1} f(u),
$$
local approximations of the time derivatives of the solution are obtained \raph{by} 
:
\begin{equation*}
u^{(k)}_{i,j} = - \mathcal{A}^{1,j}_{P}(f^{(k-1)}_{i, \bcdot}, \Delta x)
= - \frac{1}{\Delta x} \sum_{r=-P +1}^P \gamma^{1,j}_{P,r} f^{(k-1)}_{i, r}.
\end{equation*}
These approximations of the time derivatives are then adopted to compute \gio{predicted values of the flux at several time levels,  via recursive use of Taylor expansions, that will be then used to numerically approximate time derivatives of the flux at time level $n$.}
The algorithm to compute $F^p_{i+1/2}$ for the cell $i$ is then:
\begin{enumerate}
\item  {Define}
$$
f^{(0)}_{i,j}=f(u^n_{i+j}), \quad  j = -P+1, \dots, P.
$$
\item  {For $k = 2 \dots m$:}
\begin{enumerate}
\item Compute
\begin{equation*}
 u^{(k-1)}_{i,j} = - \mathcal{A}^{1,j}_{P}(f^{(k-2)}_{i,\bcdot}, \Delta x). 
\end{equation*}
\item Compute 
$$
f^{k-1,n+r}_{i,j} = f \left(  u^n_{i+j} + \sum_{l=1}^{k-1} \frac{(r \Delta t)^l}{l!} u^{(l)}_{i,j} \right), \quad  j, r = -P+1, \dots, P.
$$
\item Compute
$$
f^{(k-1)}_{i,j} =   \mathcal{A}^{k-1,0}_{P}( f^{k-1, \bcdot}_{i,j}, \Delta t),\quad  j = -P+1, \dots, P.
$$
\end{enumerate}
\item Compute $F^P_{i+1/2}$ by 
\begin{equation}
    \label{FP}
    F_{i+\ha}^P = \sum_{k=1}^{2P}\frac{\Delta t^{k-1}}{k!}f^{(k-1)}_{i+\ha},
\end{equation} where
$$ 
f^{(k-1)}_{i+\ha} = \mathcal{A}_{P}^{0,\ha}\left(f^{(k-1)}_{i,*},\Delta x\right), 
\quad 
\text{with}
\quad
\mathcal{A}_{P}^{0,\ha}\left(f^{(k-1)}_{i,*},\Delta x\right)  = \sum_{p=-P+1}^P\gamma_{P,p}^{0,\ha} \, f_{i+p}^{(k-1)}.
$$ 
\end{enumerate}

\section{CATMOOD} \label{sec:CATMOOD}


The essential idea of the MOOD technique is to apply a high-order method over the entire domain for a time step, then check locally, for each cell $i$, the behavior of the solution using some admissibility criteria such as positivity, monotonicity, physical \add{admissibility}, etc. If the solution computed in cell $i$ at time $t^{n+1}$ is in accordance with the selected criteria, it is kept. Otherwise it is \raph{locally} recomputed with a  lower order numerical method. This operation is repeated until acceptability, or when a robust first order scheme is 
\raph{employed.}

Bearing this in mind, the idea is to design a cascade of CAT methods in which the order is locally adjusted according to some \aposteriori admissibility criteria thus creating a new family of adaptive CAT methods called CAT-MOOD schemes.

\subsection{MOOD admissibility criteria} \label{ssec:MOOD}
\begin{figure}[!ht]
    \centering    
     \begin{subfigure}[b]{0.76\textwidth}
         \centering
    \includegraphics[width=0.99\textwidth]{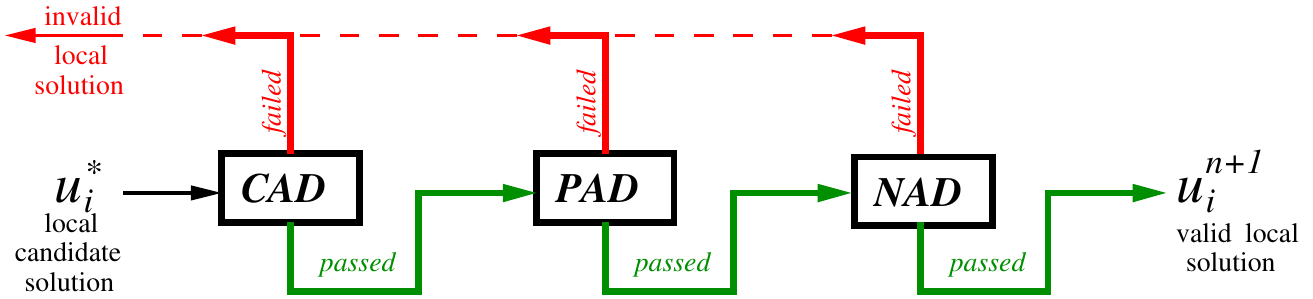}
    \caption{Criteria}
    \end{subfigure}
     \begin{subfigure}[b]{0.16\textwidth}
         \centering
    \includegraphics[width=0.99\textwidth,angle = 270]{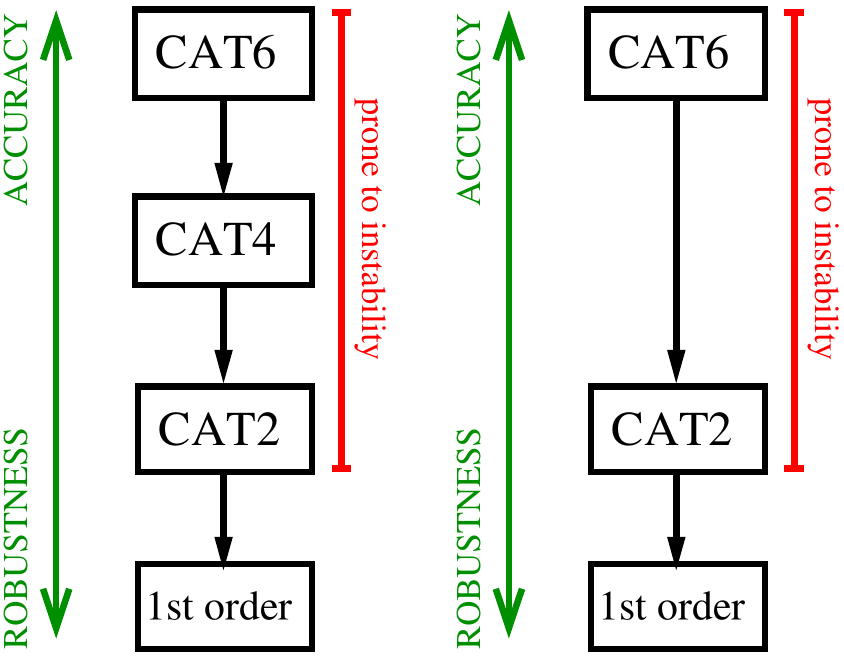}
    \caption{Cascade}
    \end{subfigure}
    \caption{Left: Detection criteria of the MOOD technique for a candidate solution $u_i^*$. \textit{Computer Admissible Detector (CAD)}, \textit{Physical Admissible Detector (PAD)}  and \textit{Numerical Admissible Detector (PAD)} ---
    Right: Order cascades of CAT schemes used in the MOOD procedure. Starting from the most accurate one, CAT6, downgrading to lower order schemes, and, at last to a $1$st order accurate scheme employed to ensure robustness.}
   \label{fig:MOOD_chain_cascade}
\end{figure}

Following \cite{CDL2,CDL3}, we select three different admissibility criteria 
 which are used to check the admissibility  of a candidate numerical solution $\left\{ u_i^{n+1} \right\}_{1\leq i \leq N}$:
\begin{enumerate}
    \item \textit{Physical Admissible Detector (PAD)}: The first detector checks the physical validity of the candidate solution. In particular, this  detector reacts to negative solution when a variable cannot take negative values: this is the case of the pressure $p$ and density $\rho$ for the 1D Euler system 
    in compliance with (\ref{eq:physical_states}). 
    \item \textit{Numerical Admissible Detector (NAD)}: This criterion is used to ensure the essentially non-oscillatory (ENO) character of the numerical solution, namely that no large and spurious minima or maxima are introduced locally in the solution. To do this, the following relaxed variant of the Discrete Maximum Principle (see \cite{Ciarlet,CDL1}) is \raph{considered}:
    $$\min_{j\in \mathcal{C}_i^P}(u_{j}^n) - \delta_i^n \le u_{i}^{n+1} \le 
    \max_{j \in \mathcal{C}_i^P}(u^n_{j}) + \delta_i^n, $$
    where {$\mathcal{C}_i^P = \{i-P,\ldots,i+P\}$} 
     is the $(2P +1)$-point centered stencil and $\delta_i$ is a parameter that avoids wrong detections in flat region. Here,  $\delta_i$ is set as:
    \begin{equation}
        \delta_i^n = \max\left(\gio{{\rm tol}_1},\gio{{\rm tol}_2},\max_{j\in \mathcal{C}_i^P}u^n_j - \min_{j\in \mathcal{C}_i^P}u^n_j\right). 
        \label{eq:delta}
    \end{equation}
 In the Euler test of Sec.~\ref{sec:numerics} the relaxed discrete maximum principle is only computed for density $\rho$ and pressure $p$. 
    \item \textit{Computer Admissible Detector (CAD)}: The last criterion detects 
    undefined or unrepresentable quantities, usually
    not-a-number \texttt{NaN} or infinity quantity \raph{which may appear}, for instance, when a division by zero \gio{is encountered}.
\end{enumerate}  
The order in which these three criteria are applied to the candidate solution is showed in Figure~\ref{fig:MOOD_chain_cascade}-left. If one of the criteria is not satisfied, the cell is marked \raph{as 'failed'}.
A new candidate solution is then computed in marked cells using a lower order method and it is checked again.

\subsection{CAT scheme with MOOD limiting} \label{ssec:CATMOOD} 
In this work, CAT methods are used as high-order methods within the MOOD strategy: the natural idea would be, \raph{given a target order $2P$}, to use the \raph{following} cascade of numerical methods 
$$\text{CAT}2P \to \text{CAT}2(P-1) \to \dots \text{CAT}2 \to \text{First-order method}$$ to  obtain a method with  order of accuracy $2P$ in smooth regions and an essentially non-oscillatory solution close to discontinuities or large gradients. 
\raph{However,}
in order to reduce the computational cost, the following cascade has been 
\raph{preferred}
in the numerical tests below
$$\text{CAT}6 \to \text{CAT}2  \to \text{First-order method}$$ 
where the first-order method is Rusanov for scalar problems and HLL for Euler equations: see Figure~\ref{fig:MOOD_chain_cascade}-right \raph{for an illustration}. Therefore, the expected order of accuracy in smooth regions is 6. The resulting scheme is referred to as CATMOOD6. 
\section{Numerical test cases}  \label{sec:numerics}

\gio{Three tests are considered, namely linear advection equation, to check the numerical order of accuracy,  Burgers equation, and Euler equations of gas dynamics, to check shock-capturing capability of the method. 
In all our tests we used the following parameters: 
${\rm tol}_1 = 10^{-4}$ and 
${\rm tol}_2 = 10^{-3}$.}

\subsection{Scalar linear equation with smooth initial condition} \label{ssec:linear}
Let us consider the linear conservation law \eqref{sec:CAT_gov_equ} with \gio{$f(u) = u$, and periodic boundary conditions in $[0,2\pi]$}. The
smooth initial condition is given by
\begin{equation}\label{smooth_1}
u(x,0) = u_0(x) = \ha\sin(x) + 1.
\end{equation} 
The goal of this test is to check and compare the empirical order of accuracy of CATMOOD6. 
For smooth solutions, the detectors are expected not to spoil the sixth-order of accuracy of the CAT6 \raph{scheme}. We run this test case in the interval $[0,2\pi],$ with final time $t_{fin} = 1$, CFL$=0.9,$ and periodic boundary conditions. We apply CAT6 and CATMOOD6 method on successive refined uniform meshes going from $10$   to  $320$ cells. As expected the empirical order of accuracy of both CAT and CATMOOD is 6: see Table \ref{tab:linear_convergence}. 
\raph{We observe that below $N=80$ cells, the mesh is not fine enough to allow for a clean limiting.
\gio{Notice that this threshold depends on the parameteters ${\rm tol}_1$ and ${\rm tol}_2$ adopted in \eqref{eq:delta}. Larger values of these parameter will lower the threshold to smaller values of $N$.} 
Moreover, in this test, CATMOOD6 is $1.5$ more expensive than CAT6.}
\begin{table}[!ht]
    \centering
    \begin{tabular}{|c||ccc|ccc|}
    \hline \multicolumn{7}{|c|}{\textbf{Linear equation - Error, Rate of convergence, CPU time}} \\
\hline   
\hline   
& \multicolumn{3}{c|}{\textbf{CAT6}} &\multicolumn{3}{c|}{\textbf{CATMOOD6}}\\
 $N$ &  $L^1$ error &  order & CPU time & $L^1$ error &  order & CPU time \\ \hline
10 & 4.27$\times 10^{-5}$  &  ---  & 0.0073    & 8.34$\times 10^{-3}$ & ---  & 0.016    \\
20 & 6.93$\times 10^{-7}$    &  5.95  & 0.020    & 3.18$\times 10^{-3}$ & 1.39  & 0.031    \\
40 & 9.64$\times 10^{-9}$    &  6.17  & 0.038     & 3.16$\times 10^{-4}$ & 3.33  & 0.053    \\
80 & 1.43$\times 10^{-10}$   &  6.07  & 0.068    &  2.48$\times 10^{-6}$ &  3.67  & 0.094    \\
160 & 2.09$\times 10^{-12}$  &  6.09  & 0.13    & 2.09$\times 10^{-12}$ & 23.50  & 0.2    \\
320 & 3.29$\times 10^{-14}$  &  5.99  & 0.25    & 3.31$\times 10^{-14}$ & 5.99   & 0.38    \\
& \text{Expected} & 6 &
& \text{Expected} & 6 &
\\
    \hline 
    \end{tabular}
    \caption{Linear scalar equation \ref{ssec:linear}.  $L^1-$norm errors between the numerical solution and the exact solution of the linear equation at $t_{\text{final}} = 1$ on uniform Cartesian mesh and CFL$=0.9.$
   }
    \label{tab:linear_convergence}
\end{table}

\vspace{-0.6cm}

\subsection{Burgers' equation with non-smooth initial condition} \label{ssec:burgers}
\begin{figure}[!ht]
\hspace{-1.cm}
    \centering    
    \begin{subfigure}[b]{0.48\textwidth}
         \centering
    \includegraphics[width=1.12\textwidth]{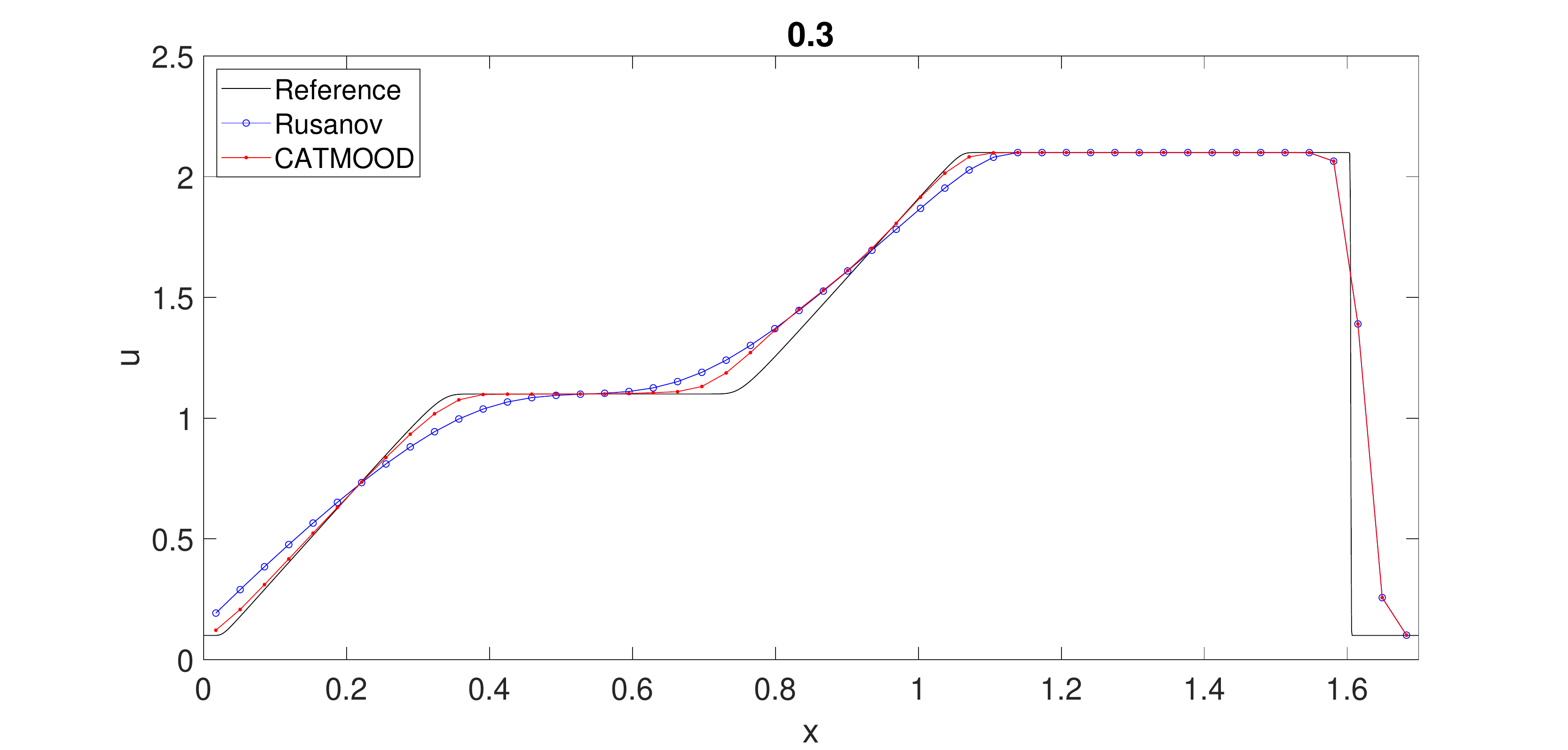}
    \vspace{-0.5cm}
    \caption{Numerical solutions at time $t=0.3$.}
    \end{subfigure}
     \begin{subfigure}[b]{0.48\textwidth}
         \centering
    \includegraphics[width=1.12\textwidth]{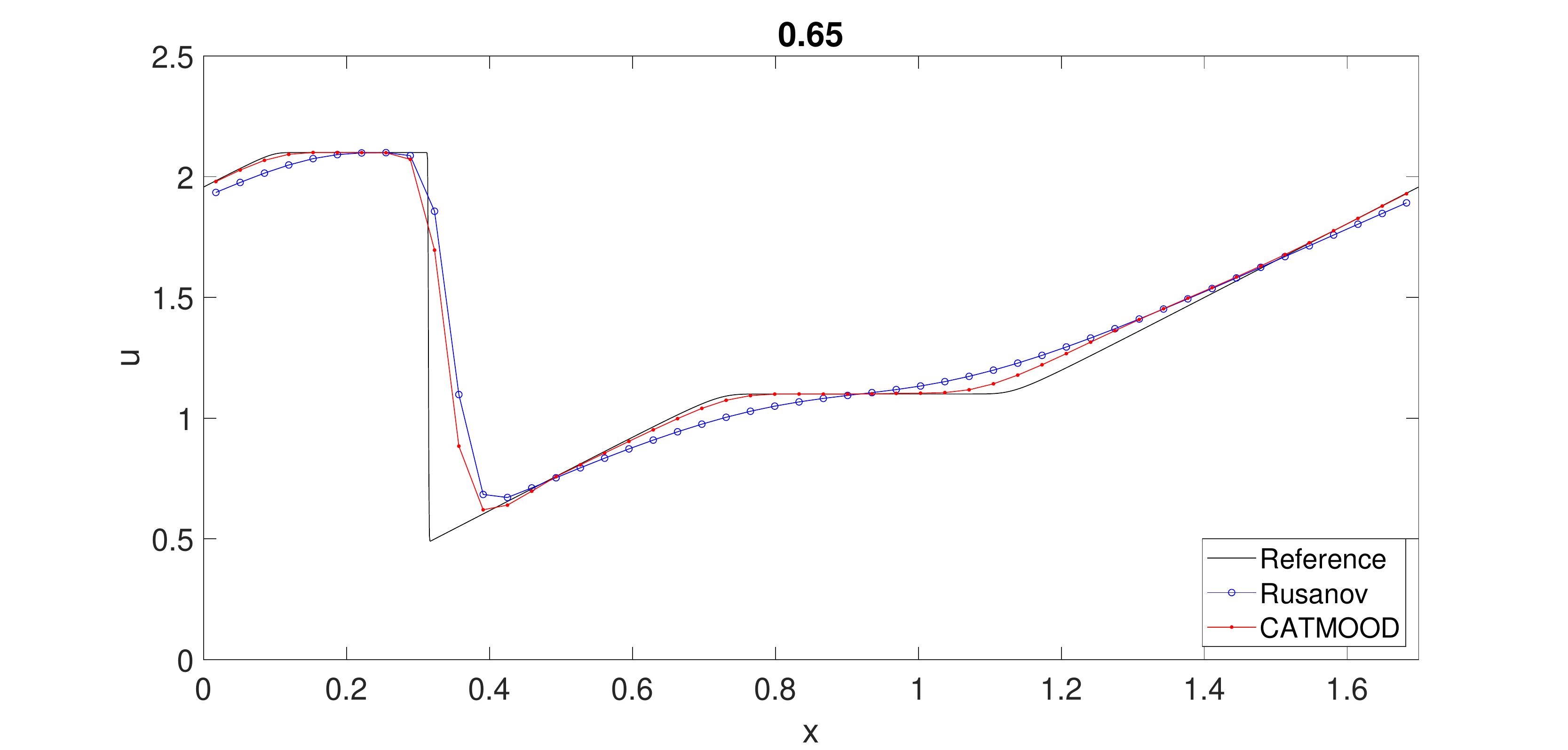}
    \vspace{-0.5cm}
    \caption{Numerical solutions at time $t=0.65$.}
    \end{subfigure}
    \caption{Burgers' equation with non-smooth initial condition \ref{ssec:burgers} --- Left: Numerical solutions at time $t = 0.3$ --- Right: Numerical solutions at final time $t_{fin} = 0.65$ obtained with Rusanov flux and CATMOOD6 on $50$ uniform mesh and CFL$=0.9$. Rusanov flux scheme has been used as first order. The reference solution has been obtained with the 1st order Rusanov scheme on $2000$ uniform cells. }
   \label{fig:burgers}
\end{figure}

Let us consider \eqref{sec:CAT_gov_equ} with \gio{$f(u) = u^2/2$, periodic boundary conditions,} and non-smooth initial condition
\begin{equation} \label{square_step_test}
u_0(x) = \begin{cases}
1.1 \quad \;\;\;\mathrm{if}\quad 0\le x\le \ha;\\
2.1 \quad\;\; \;\mathrm{if}\quad \ha< x < \frac{3}{2};\\
0.1 \quad \;\;\; \mathrm{if} \quad \frac{3}{2}\le x\le\frac{17}{10}. 
\end{cases}
\end{equation} 
We run this test case with first-order Rusanov-flux and  CATMOOD6 method in the interval $[0,1.7],$
with final time $t_{fin} = 0.65$, 
using a 50-cell mesh, CFL$=0.9,$ and periodic boundary conditions.
Figure~\ref{fig:burgers} shows the initial condition and the numerical solutions obtained with both methods: it can be seen that CATMOOD6 provides a non-oscillatory solution \gio{still providing better resolution than the first order method. Notice that larger values of the thresholds tol$_1$ and tol$_2$ will decrease dissipation, but may not be sufficient to avoid creation of spurious oscillations.}

\vspace{-0.6cm}

\subsection{Euler system: Sod problem} \label{ssec:Sod}
Let us consider Euler equations \eqref{sec:CAT_gov_equ} with SOD initial condition 

\begin{equation} \label{Sod_IC}
	(\rho, \raph{u}, p)
	= \left\{
	\begin{array}{ll}
	\displaystyle (1,0,1) & \mbox {if } x \le 0, \\
	\displaystyle (0.125,0,0.1) & \mbox {if } x > 0.
	\end{array}\right.
\end{equation}

\begin{figure}[!ht]
\hspace{-1.cm}
    \centering    
    \begin{subfigure}[b]{0.48\textwidth}
         \centering
    \includegraphics[width=1.12\textwidth]{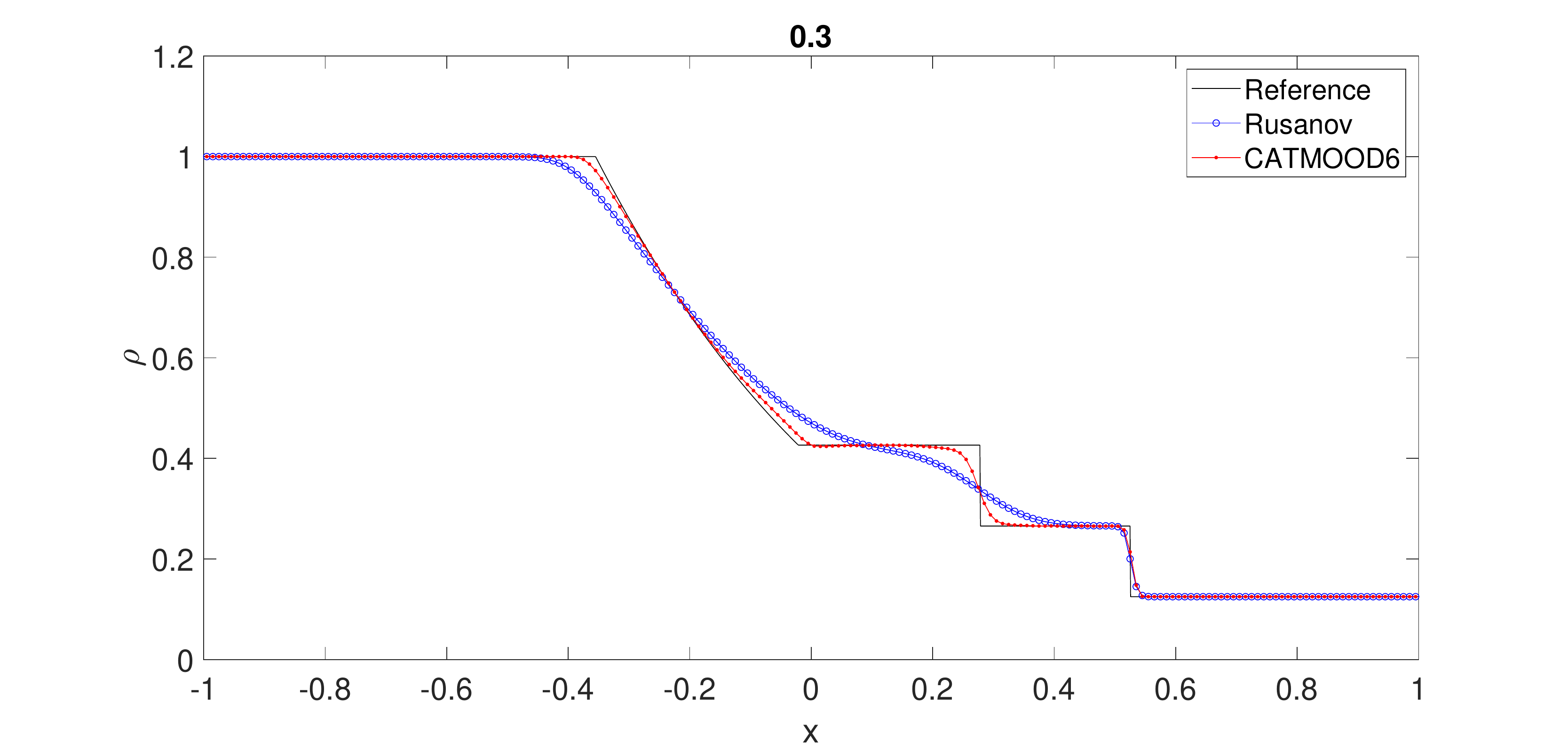} 
    \vspace{-0.5cm}
    \caption{Density.}
    \end{subfigure}
     \begin{subfigure}[b]{0.48\textwidth}
         \centering
    \includegraphics[width=1.12\textwidth]{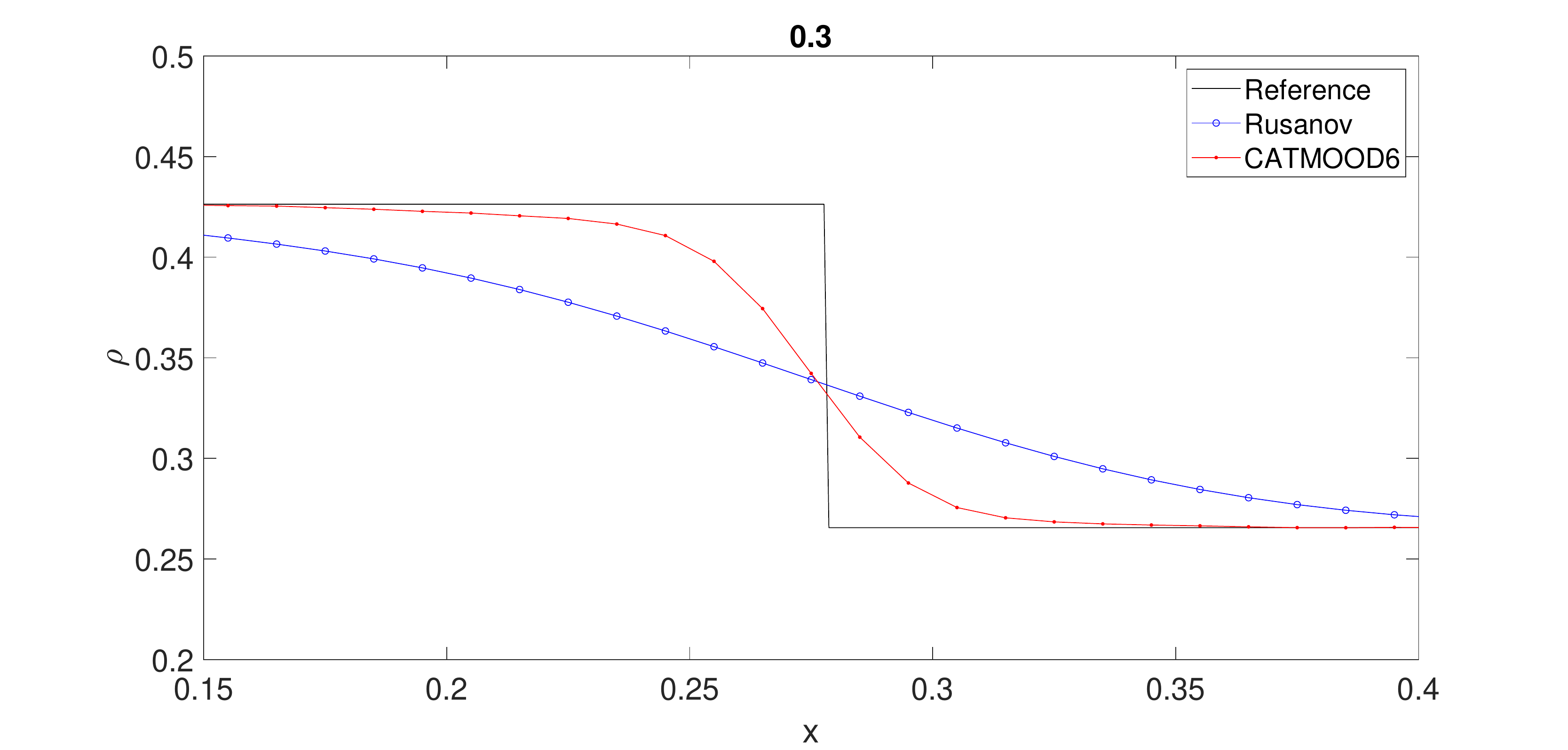}
    \vspace{-0.5cm}
    \caption{Zoom of the density on shock}
    \end{subfigure} 
    \\
    \hspace{-1.cm}
    \centering    
    \begin{subfigure}[b]{0.48\textwidth}
         \centering
    \includegraphics[width=1.12\textwidth]{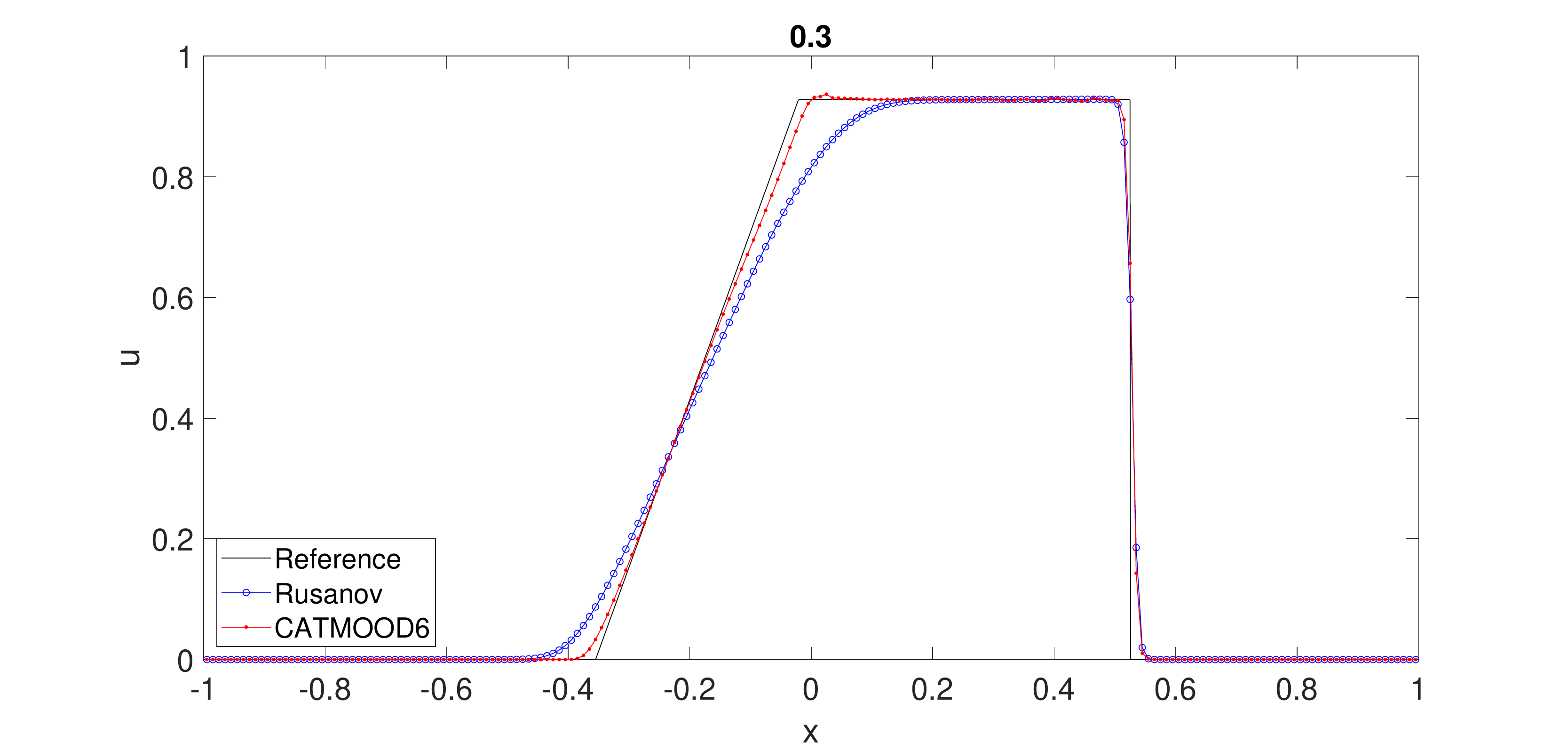} 
    \vspace{-0.5cm}
    \caption{Velocity.}
    \end{subfigure}
     \begin{subfigure}[b]{0.48\textwidth}
         \centering
    \includegraphics[width=1.12\textwidth]{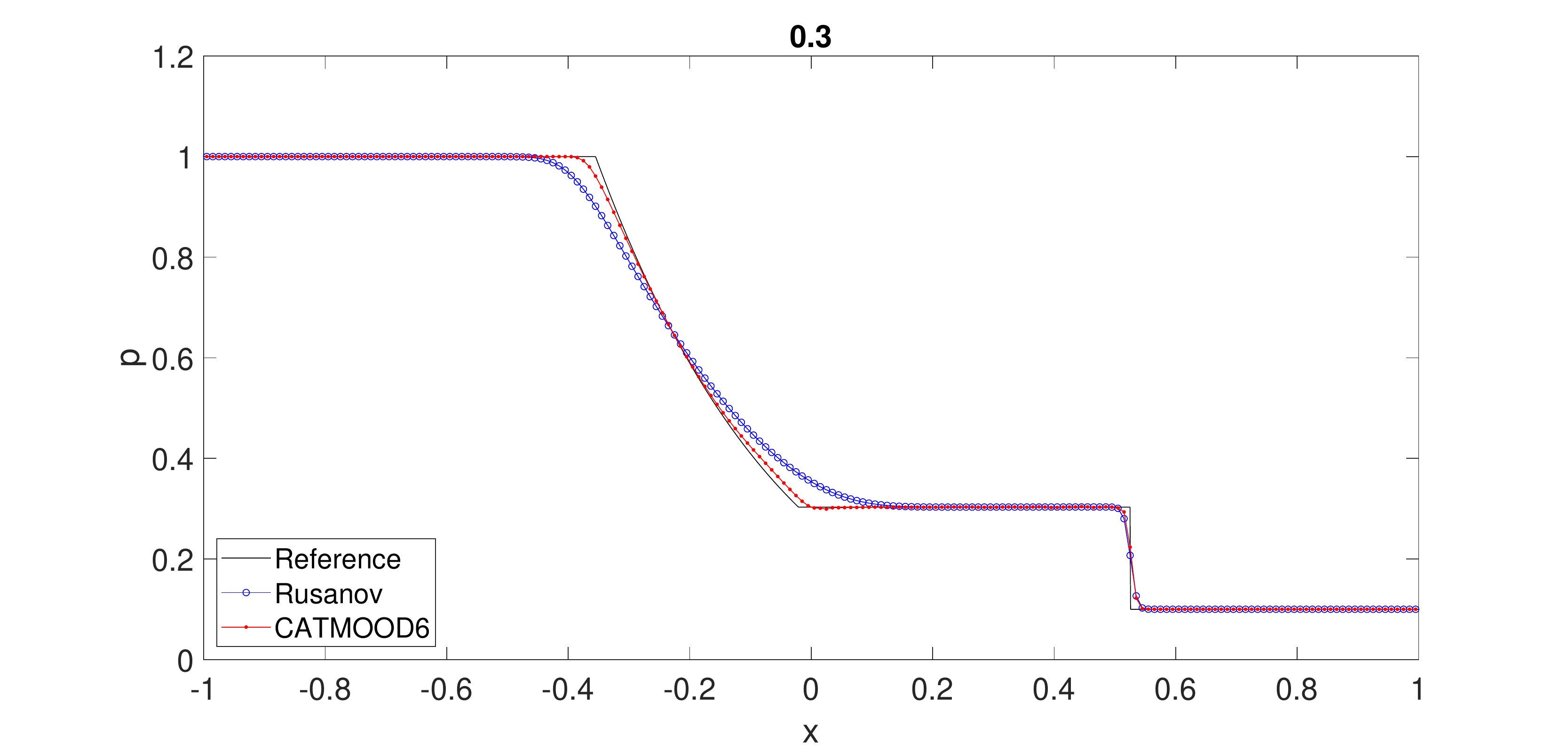}
    \vspace{-0.5cm}
    \caption{Pressure}
    \end{subfigure}
    \caption{Euler system: Sod problem \ref{ssec:Sod} --- Left-up: Numerical solutions for density $\rho$ at final time $t=0.3$ obtained with \raph{1st order HLL} 
    flux and CATMOOD6, on $200$ uniform mesh and CFL$=0.9$. --- Right-up: Zoom of the numerical solutions and reference one. --- 
    Left/right-down: same as left-up but for velocity and pressure variables. The reference solution is the exact one.}
   \label{fig:Sod}
\end{figure}
We run this test case with first-order Rusanov-flux and CATMOOD6 method in the interval $[-1,1],$ with final time $t_{fin} = 0.3$. We use  a $200-$cell mesh, CFL$=0.9,$ and free boundary conditions.
Figure~\ref{fig:Sod} shows the numerical and  the exact solutions for density, velocity and pressure obtained with both methods and a zoom of the density close to the shock. It can be seen that, as expected, CATMOOD6 solution shows a better resolution of rarefaction and contact waves.

\section{Conclusion and perspectives} \label{sec:conclusion}
In this paper we have presented a combination between the \aposteriori shock-capturing MOOD technique  and the one-step high-order finite-difference CAT2P schemes. CAT2P schemes are of order $2P$ on smooth solutions, but, \raph{an} extra dissipative mechanism must be supplemented to deal with steep gradients or discontinuous solutions.
In this work we rely on an \aposteriori MOOD paradigm which computes an unlimited high-order candidate solution at time $t^{n+1}$, and, further detects troubled cells which are recomputed with a lower-order scheme throughout a family of detectors.
For a proof of concept, we tested the so-called \raph{CATMOOD6 scheme based on the } 'cascade': 
CAT6$\rightarrow$CAT2$\rightarrow \,$ 1st, where the last scheme is a first order robust scheme.
\raph{This scheme has been challenged}
on a test suite of smooth solutions (linear scalar equation), simple shock waves (Burgers' equation), and complex self-similar solutions involving contact, shock and rarefaction waves (Sod problem).
In all test cases, CATMOOD6 has preserved the accuracy on smooth parts of the solutions, an essentially-non-oscillatory behavior close to steep gradients, and, \raph{always produces} a physically valid solution.
\gio{A detailed analysis of the improvement in the efficiency over standard slope limiters has not been performed. 
A two dimensional implementation, as well as other generalizations and  improvements are  under way.}

\vspace{-0.75cm}
\bibliographystyle{plain}
\bibliography{biblio}

\end{document}